\documentclass[10pt]{article}
\usepackage{amsmath,amssymb,amsthm,amscd}
\numberwithin{equation}{section} 

\def\p{\partial}

\def\i{\sqrt{-1}}

\def\cC{{\cal C}}

\def\cA{{\mathcal A}}
\def\cC{{\mathcal C}}

\def\cF{{\mathcal F}}

\newtheorem{prop}{Proposition}[section]
\newtheorem{theo}[prop]{Theorem}
\newtheorem{lemma}[prop]{Lemma}

\newtheorem{rem}[prop]{Remark}

\let\lra=\longrightarrow

\def\mapright\#1{\,\smash{\mathop{\lra}\limits^{\#1}}\,}

\begin{document}
\bibliographystyle{plain}
\title{Remarks on the existence of bilaterally symmetric extremal K\"ahler
metrics on $\mathbb{CP}^2\sharp 2\overline{ \mathbb{CP}^2}$}
\author{He, Weiyong}
\date{May 28th, 2007}
\maketitle

The study of extremal K\"ahler metric is initiated by the seminal
works of Calabi \cite{Calabi01}, \cite{Calabi02}.  Let $(M,
[\omega])$ be a compact K\"ahler manifold with fixed K\"ahler
class $[\omega]$. For any K\"ahler metrics $g$
in the fixed K\"ahler class $[\omega]$, the Calabi energy $\cC(g)$  is defined as
\[\cC(g)=\int_Ms^2d\mu,\] where $s$ is the scalar curvature of
$g$.  The extremal K\"ahler metric is the critical
point of the Calabi energy.  The Euler-Lagrange equation is
\[\bar\p \nabla^{1, 0} s=0.
\]
In other words, $\Xi = \nabla ^{1, 0}s$ is a holomorphic vector
field  (we call it extremal vector field from now on). From PDE
point of view, the existence of the extremal metric is to solve a
6th order nonlinear elliptic equation.   According to Chen
\cite{Chen} (c.f. Donaldson \cite{Donaldson} for algebraic case),
there is {\it a priori} greatest lower bound for the Calabi energy
in any fixed K\"ahler class.   This {\it a priori} lower bound can be computed explicitly  as
\[
\cA([\omega])=\frac{(c_1\cdot
[\omega])^2}{[\omega]^2}-\frac{1}{32\pi^2}\cF(\Xi, [\omega]),
\] where $\cF(\Xi, [\omega])$ is the Futaki invariant of class $  [\omega]$.  Note that the extremal vector field
$\Xi$ is determined \cite{Futaki-Mabuchi} up to conjugation
without the assumption of the existence of an extremal metric.

By E. Calabi \cite{Calabi02}, extremal K\"ahler metrics minimizes the Calabi energy locally.   By X.X. Chen (\cite{Chen}) and S.K. Donaldson (\cite{Donaldson}),  we know
\[ \cA([\omega])\leq \frac{1}{32\pi^2}\min_{g\in
[\omega]}\cC(g),\] where the equality holds when there is an
extremal K\"ahler metric in $
[\omega].$\\

In an amazingly  beautiful work, Chen-LeBrun-Weber
\cite{chen-lebrun-weber} proved the existence of bilaterally
symmetric extremal K\"ahler  metrics on $\mathbb{CP}^2\sharp
2\overline{\mathbb{CP}^2}$ by global deformation method.  More
strikingly, it contains an extremal class where the extremal
metric is conformal to an Einstein metric with positive  scalar
curvature. $\mathbb{CP}^2\sharp 2\overline{\mathbb{CP}^2}$ can  be
also described as $\mathbb{CP}^1\times \mathbb{CP}^1$ blowing up
at one point. We use $F_1, F_2$ to denote the Poinc\'are dual of
two factors $\mathbb{CP}^1$ in $\mathbb{CP}^1\times
\mathbb{CP}^1,$ and $E$ denotes the exceptional divisor. The term
``bilaterally symmetric" is introduced in \cite{chen-lebrun-weber}
to describe the K\"ahler class which are invariant under the
interchange $F_1\leftrightarrow F_2$. The ``bilaterally symmetric"
class can be described by $[\omega]_x=(1+x)(F_1+F_2)-xE$ for
$0<x<\infty.$ Let $f(x)=\cA([\omega]_x)$, and it is shown that
$f(x)<9$ (c.f. \cite{chen-lebrun-weber}). Set $L$ to be the
smallest number of $f^{-1}(8)$, Chen-LeBrun-Weber
\cite{chen-lebrun-weber} proved the following theorem regarding
the existence of extremal K\"ahler
metrics\\

\noindent {\bf Theorem A} \cite{chen-lebrun-weber} {\it For any $x
\in (0, L)$, let $[\omega]_x=(1+x)(F_1+F_2)-xE$ denote the
K\"ahler class of on $M=\mathbb{CP}^2\sharp 2\overline{
\mathbb{CP}^2}$, then there is an extremal metric in $[\omega]_x$
for any $x\in (0, L).$}

Their method is through large scale  deformation. The existence of
extremal K\"ahler metrics is promised by the results of
Arezzo-Pacard-Singer \cite{APS} when $x$ is small enough (also for
$x$ big enough).  According to LeBrun-Simanca
\cite{LeBrun-Simanca}, the set which admits extremal K\"ahler
metric is open. Following the work of  Chen-Weber
\cite{Chen-Weber} on moduli space of  extremal K\"ahler metrics 
in complex surface, a sequence of bilaterally symmetric extremal
metrics will converge to an extremal metric with finite orbifold
points. However the orbifold singularities can only arise as a
very specific mechanism of curvature concentration for critical
metrics \cite{Anderson}, \cite{Tian-Viaclovsky},
\cite{Chen-Weber}. The key idea of Chen-LeBrun-Weber
\cite{chen-lebrun-weber} is thorough careful analysis of the
bubble formation and they conclude that, for bubble to arise, the
original K\"ahler class must admit some Lagrange cycle with
negative self-intersection number. And they show that when
$f(x)<8$, there is no such Lagrange cycle. It follows that the
orbifold singularities will never occur.\\

Inspired by the idea of \cite{chen-lebrun-weber}, we extend their
result to show that the existence of bilaterally symmetric
extremal K\"ahler metrics on $\mathbb{CP}^2\sharp 2\overline{
\mathbb{CP}^2}$ for any $x\in (0, \infty)$ in this short note. The
readers are enthusatically referred to \cite{chen-lebrun-weber}
for the
 historic background of this problem as well as an
excellent list of references.  Following the scheme in
(\cite{chen-lebrun-weber}),  we show that
\begin{theo} For any $x\in (0, \infty)$, let
$[\omega]_x=(1+x)(F_1+F_2)-xE$ denote the K\"ahler class of on
$M=\mathbb{CP}^2\sharp 2\overline{ \mathbb{CP}^2}$, then there is
an extremal metric in $[\omega]_x$ for any $x\in (0, \infty).$
\end{theo}

We keep the notations of \cite{chen-lebrun-weber}.  Our observation is that,  without assuming
 $\cA([\omega]_x) < 8$, the proposition
(\cite{chen-lebrun-weber}, Proposition 26) still holds.
\begin{prop}Let $g_i$ be a sequence of unit-volume bilaterally
symmetric extremal K\"ahler metrics on $(M, J)=\mathbb{CP}^2\sharp
2\overline{\mathbb{CP}^2}$ such that the corresponding K\"ahler
class
\[[\omega_i]=\frac{(1+x_i)(F_1+F_2)-x_iE}{\sqrt{1+2x_i+x_i^2/2}}\]
satisfy $A\leq x_i\leq B$, where $A<B$ are any two fixed positive
number. Then there is a subsequence $g_{i_j}$ of metrics and a
sequence of diffeomorphisms $\Psi_j: M\rightarrow M$ such that
$\Psi_j^{*}g_{i_j}$ converges in the smooth topology to an
extremal K\"ahler metric on the smooth 4-manifold $M$ compatible
with some complex structure $\tilde{J}=\lim_{j\rightarrow
\infty}\Psi_{j*}J.$
\end{prop}
Recall for a compact smooth 4-manifold $(M, g)$ the Gauss-Bonnet
formula says
\[
\frac{1}{8\pi^2}\int_M\left(|W_{+}|^2+|W_{-}|^2+\frac{s^2}{24}-\frac{|
Ric_0|^2}{2}\right)d\mu=\chi(M)
\]
and the signature formula reads
\[
\frac{1}{12\pi^2}\int_M(|W_{+}|^2-|W_{-}|^2)d\mu=\tau(M).\] If
$(X, g_\infty)$ is any ALE 4-manifold with finite group
$\Gamma\subset SO(4)$ at infinity, then the Gauss-Bonnet formula
becomes
\[
\frac{1}{8\pi^2}\int_X\left(|W_{+}|^2+|W_{-}|^2+\frac{s^2}{24}-\frac{|
Ric_0|^2}{2}\right)d\mu_{g_\infty}=\chi(X)-\frac{1}{|\Gamma|}
\]
and the signature formula becomes
\[
\frac{1}{12\pi^2}\int_X(|W_{+}|^2-|W_{-}|^2)d\mu_{g_\infty}=\tau(X)+
\eta(S^3/\Gamma),
\]
where $\chi(X)$ is the Euler characteristic of non-compact
manifold $X$ and $\eta(S^3/\Gamma)$ is called $\eta$ invariant.
When $(X, g_\infty)$ is scalar flat K\"ahler, the formulas
simplify to
\[
\frac{1}{8\pi^2}\int_X\left(|W_{-}|^2-\frac{|Ric_0|^2}{2}\right)d\mu_
{g_\infty}=\chi(X)-\frac{1}{|\Gamma|}
\]
and

\[
-\frac{1}{12\pi^2}\int_X|W_{-}|^2d\mu_{g_\infty}=\tau(X)+\eta(S^3/
\Gamma).
\]

Our first observation is that the lemmas
(\cite{chen-lebrun-weber}, Lemma 21 and Lemma 22) hold without the
assumption on $\cA([\omega]).$
\begin{lemma}$(X, g_\infty)$ is the deepest bubble. Then $X$ is
diffeomorphic to a region of $M$ which is invariant under
$F_1\leftrightarrow F_2$, and this $\mathbb{Z}_2$ action induces a
holomorphic isometric involution of $(X, g_\infty).$
\end{lemma}
\begin{proof}
By the signature formula, we have that
\[
\int_M|W_{-}|^2d\mu=-12\tau(M)+\int_M|W_{+}|^2d\mu=12\pi^2+\int_M\frac
{s^2}{24}d\mu
\]
for any K\"ahler metrics on $M=\mathbb{CP}^2\sharp
2\overline{\mathbb{CP}^2}$. For any bilaterally symmetric K\"ahler
class $[\omega]$ on $\mathbb{CP}^2\sharp
2\overline{\mathbb{CP}^2}$, $\cA([\omega])<9$. Thus any
bilaterally symmetric extremal K\"ahler metrics satisfy
\[
\int_M|W_{-}|^2d\mu<12\pi^2+\frac{9}{24}32\pi^2=24\pi^2.
\]

When $|\Gamma|\geq 2$, since $b_1(X)=b_3(X)=0$ and $b_2(X)>0$.
Hence $\chi(X)\geq 2$, and the Gauss-Bonnet formula gives that
\[
\int_X|W_{-}|^2d\mu_{g_\infty}\geq 8\pi^2(2-1/2)\geq 12\pi^2.
\]
When $|\Gamma|=1$, the signature formula gives that
\[
\int_X|W_{-}|^2d\mu =12\pi^2.
\]
And then the same argument of (\cite{chen-lebrun-weber} Lemma 21)
applies.
\end{proof}
\begin{lemma}Let $(X, g_\infty)$ be the deepest bubble. If
$b_2(X)=1$, then $X$ must be diffeomorphic to the line bundle of
degree $-k$ over $\mathbb{CP}^1$ for $1\leq k\leq 5$.
\end{lemma}
\begin{proof}The proof follows ({\cite{chen-lebrun-weber}}, Lemma   23). Since
$X$ is diffeomorphic to the line bundle of degree $-k$ over
$\mathbb{CP}^1$ for some $k>0$. If $C$ denotes the homology class
of the zero section, the Poincar\'e dual of $c_1$ is the rational
homology class $\frac{k-2}{k}C$ and it follows that
\begin{equation*}
\int_{X}|Ric_0|^2d\mu_{g_\infty}=-8\pi^2c_1^2=8\pi^2\frac{(k-2)^2}{k}.
\end{equation*}
Any bilaterally symmetric extremal K\"ahler metrics satisfy
\[
\int_M|Ric_0|^2d\mu=\frac{1}{4}\int_Ms^2d\mu-8\pi^2c_1^2(M)<16\pi^2.
\]
It follows that $k\leq 5.$ \end{proof} (\cite{chen-lebrun-weber}
Lemma 22) holds also.
\begin{lemma}Let $(X, g_\infty)$ be the deepest bubble. If
$b_2(X)=2$,then $\Gamma\cong \mathbb{Z}_3$, and $X$ has
intersection form
\[\left(
\begin{array}{cl}
-2  \hspace{2mm}&1\\
1   \hspace{2mm}&-2\\
\end{array}
\right).\]
\end{lemma}
\begin{proof}
Since $b_2(M_\infty)=2$, the Gauss-Bonnet and signature formula
give that
\[
\frac{1}{12\pi^2}\int_{X}|W_{-}|^2d\mu_{g_\infty}=2-\eta(S^2/\Gamma)
\]
and
\[
\frac{1}{8\pi^2}\int_{X}\left(|W_{-}|^2-\frac{|Ric_0|^2}{2}\right)d
\mu_{g_\infty}=3-\frac{1}{|\Gamma|}.
\]
It follows that
\begin{equation}\label{2-1}
\frac{3}{2}\eta(S^3/\Gamma)+\frac{1}{16\pi^2}\int_{X}|Ric_0|^2d\mu_{g_
\infty}=\frac{1}{|\Gamma|}.
\end{equation}
And we know that
\[
\int_M|W_{-}|^2d\mu=12\pi^2+\int_M\frac{s^2}{24}d\mu<24\pi^2,
\]
it follows that
\[
\frac{1}{16\pi^2}\int_X\frac{|Ric_0|^2}{2}<\frac{1}{|\Gamma|}
\]
and
\[
\eta(S^3/\Gamma)>0.
\]
Since Lemma 0.3 shows that we still have a $\mathbb{Z}_2$ action
which interchanges the two totally geodesic $\mathbb{CP}^1s$ which
generate $H^2(X, \mathbb{Z})$. The argument in
(\cite{chen-lebrun-weber} Lemma 22) applies and so the
intersection form of $X$ must be given by \[\left(
\begin{array}{cl}
-k  \hspace{2mm}&1\\
1   \hspace{2mm}&-k\\
\end{array}
\right)\]for some $k\geq 2$ and $\Gamma\cong\mathbb{Z}_{k^2-1}.$
And at infinity the 3-manifold is a Lens space $L(k^2-1, k)$. In
particular $\Gamma\neq \{1\}$. Since $|\Gamma|\neq 1$,  by
(\ref{2-1}) we  get that
\[
\eta(S^3/\Gamma)\leq \frac{1}{3}. \]It means that
\[
0<\eta(S^3/\Gamma)\leq \frac{1}{3}.
\]
For the Lens space $L(k^2-1, k)=S^3/\Gamma$, the $\eta$-invariant
is given by \cite{AtiyahPS},
\begin{eqnarray}\label{2-2}
\eta(S^3/\Gamma)&=&-\frac{1}{|\Gamma|}\sum_{i=1}^{k^2-2}\cot{\frac{i
\pi}{k^2-1}}\cot{\frac{ki\pi}{k^2-1}}\nonumber\\
&=&-\frac{1}{k^2-1}{\left(\frac{2}{3}k^3-2k^2+2\right).}
\end{eqnarray}

It follows that $k=2$ and $\Gamma\cong \mathbb{Z}_3$.
\end{proof}
\begin{rem}
In this case, one can calculate the first Chern class in stead of
the $\eta$-invariant as in Lemma 0.4. And the Poincar\'e dual of
the  first Chern class is the rational homology class
\[
\frac{k-2}{k-1}(E_1+E_2),
\]
where $E_1, E_2$ are two totally geodesic $\mathbb{CP}^1$ and they
have intersection form\[\left(
\begin{array}{cl}
-k  \hspace{2mm}&1\\
1   \hspace{2mm}&-k\\
\end{array}
\right).\] But the calculation of the eta-invariant will have
independent interest for lens spaces. The formula is given by
\cite{AtiyahPS}.  We carry out the example for lens spaces
$L(k^2-1, k)$.
\end{rem}

\begin{lemma}Under the assumption of Proposition 0.2, for any $A, B$   fixed,
$X$ can not be as in Lemma 0.4 and Lemma 0.5.
\end{lemma}
\begin{proof}The proof follows exactly (\cite{chen-lebrun-weber} Lemma
25). Since the limit metric $g_\infty$ on $X$ is by construction a
pointed limit of larger and larger rescalings of the metrics
$g_i$, the generators of $H_2(X, \mathbb{Z})$ must arise from
smooth 2-sphere $S_i\subset M$ whose areas with respect to $g_i$
tend to zero as $i\rightarrow \infty.$ When $b_2(X)=1$, let $S_i$
be the smooth 2-sphere corresponding to the zero section
$\mathbb{CP}^1$; when $b_2(X)=2$, let $S_i$ be a 2-sphere
corresponding to one of the two $\mathbb{CP}^1$ generators, and
$\tilde{S_i}$ is the reflection under $F_1\leftrightarrow F_2.$
Take $\Sigma= [S_i]\in H^2(M, \mathbb{Z})$ when $b_2(X)=1$, and
$\Sigma=[S_i]+[\tilde S_i]$ when $b_2(X)=2.$ Since the homology
class is $\mathbb{Z}_2$ invariant, we have
\[
[S]=m(F_1+F_2)+nE
\] for some integers $m$ and $n$ and the self-intersection
condition gave that
\[
2m^2-n^2=-k
\]
for $k\leq 5$. Now any of unit-volume bilaterally symmetric
K\"ahler classes $[\omega_i]$ is of the form
\[
[\omega_i]=\frac{(1+x_i)(F_1+F_2)-x_iE}{\sqrt{1+2x_i+x_i^2/2}},
\]
where $A\leq x_i\leq B$.

Also we know that the area of $S_i$ measured by $g_i$ goes to zero
when $i\rightarrow \infty$. By Wirtinger's inequality we can get
\[\left|[\omega_i][\Sigma]\right|<2\; {\rm area}(S_i)\rightarrow 0.\]
It follows that
\[
\frac{2m(1+x_i)+nx_i}{\sqrt{1+2x_i+x_i^2/2}}\rightarrow 0.
\]
Denote
\[
\frac{2m(1+x_i)+nx_i}{\sqrt{1+2x_i+x_i^2/2}}=\varepsilon_i,
\]
we can get that
\[
n=-2m\frac{1+x_i}{x_i}+\varepsilon_i\frac{\sqrt{1+2x_i+x_i^2/2}}{x_i}.
\]
Since $\frac{\sqrt{1+2x_i+x_i^2/2}}{x_i}$ is uniformly bounded for
$x_i\in [A, B]$ and $\epsilon_i\rightarrow 0$ when $i\rightarrow
0,$ then
\[
4m^2\left(\frac{(1+x_i)^2}{x_i^2}-\epsilon \right)-\varepsilon_iC
(\epsilon, A, B)\leq n^2\leq
4m^2\left(\frac{(1+x_i)^2}{x_i^2}+\epsilon\right)+\varepsilon_iC
(\epsilon, A, B)
\]
where $\epsilon$ is arbitrary small positive number and
$C(\epsilon)$  is independent of $i$. We can take
$\epsilon=\frac{1}{100}$ and when $i$ big enough, $C(1/100, A,
B)\varepsilon_i<1/100$, then it gives that
\[
4m^2\frac{(1+x)^2}{x^2}-2m^2-\frac{m^2}{100}\leq k+1/100.
\]
It follows that \[ (2-1/100)m^2<k+1/100.
\]
Since $k\leq 5$, it gives that $m=0, \pm1.$ But $m=0$ gives that
$n=0 $, contradiction. If $m=1,$ then $k=2, n=-2$. And $m=-1,$
then $k=2, n=2$. For any cases,
\[
|\varepsilon_i|=\left|\frac{2m(1+x_i)+nx_i}{\sqrt{1+2x_i+x_i^2/2}}
\right|=\frac{2}{\sqrt{1+2x_i+x_i^2/2}}
\]
is uniformly bounded for $x_i\in [A, B]$. Contradiction.
\end{proof}
Deepest bubbles can therefore never arise, Proposition 0.2
follows. By using the result (\cite{chen-lebrun-weber}, Theorem
27), Proposition 0.3 implies that the existence of bilaterally
symmetric extremal K\"ahler metrics in the bilaterally symmetric
K\"ahler class for any $x\in [A, B]$.
\section{Appendix}
Here we prove the identity in  (0.2.)\begin{equation}\label{1-1}
\sum_{i=1}^{k^2-2}\cot{\frac{i\pi}{k^2-1}}\cot{\frac{ki\pi}{k^2-1}}=
\frac{2}{3}k^3-2k^2+2.\end{equation} and it follows that the
eta-invariant for lens space $L(k^2-1, k)$ is
\[
-\frac{1}{k^2-1}{\left(\frac{2}{3}k^3-2k^2+2\right)}.
\]
\begin{lemma}
$k\in \mathbb{N},$
\begin{equation}{\label{1-2}}
\sin{(k+1)x}=2^k\prod_{i=0}^k\sin\left(x+\frac{i\pi}{k+1}\right).
\end{equation}
\begin{proof}
\[
2\sin x=\i (e^{-\i x}-e^{\i x})=\i e^{-\i x}(1-e^{2\i x}).
\]
It follows that
\begin{eqnarray*}
\prod_{i=0}^{k}2^{k+1}\sin
\left(x+\frac{i\pi}{k+1}\right)&=&\prod_{i=0}^{k}\left\{\i e^{-\i
(x+\frac{i\pi}{k+1})}(1-e^{2\i (x+\frac{i\pi}{k+1})})\right\}\\
&=&(\i)^{k+1}e^{-(k+1)\i x-{k\i\pi\over 2}}e^{2(k+1)\i
x}\prod_{i=0}^{k}(e^{-2\i x}-e^{\frac{2i\i \pi}{k+1}})\\
&=&\i(e^{-(k+1)\i x}-e^{(k+1)\i x})\\
&=&2\sin (k+1)x.
\end{eqnarray*}
\end{proof}
\end{lemma}
\begin{lemma}\[(k+1)\cot{(k+1)x}=\sum_{i=0}^{k}\cot{\left(x+\frac{i
\pi}{k+1}\right)}.\]
\end{lemma}
\begin{proof}
Taking derivative on both sides of (\ref{1-2}), it gives that
\[
(k+1)\cos{(k+1)x}=2^k\cos{\left(x+\frac{j\pi}{k+1}\right)}\prod_{i\neq
j}\sin{\left(x+\frac{i\pi}{k+1}\right).}
\]
Then divided by (\ref{1-1}), we get
\[(k+1)\cot{(k+1)x}=\sum_{i=0}^{k}\cot{\left(x+\frac{i\pi}{k+1}  \right)}.\]
\end{proof}
\begin{lemma}
\[
(k+1)^2\cot^2{(k+1)x}+(k+1)k=\sum_{i=0}^{k}\cot^2{\left(x+\frac{i\pi}
{k+1}\right)}.
\]
\end{lemma}
\begin{proof}In Lemma 1.2., taking derivative on both sides.
\end{proof}
\begin{lemma}\[\sum_{i=1}^{k}\cot^2{\frac{i\pi}{k+1}}=\frac{k(k-1)}{3}.
\]
\end{lemma}
\begin{proof}
In Lemma 1.3., by taking limit for $x\rightarrow 0$.
\end{proof}
Now we can prove (\ref{1-1}).
\begin{proof}When $i$ is not the multiple of $k-1$, then
\begin{eqnarray*}
\cot{\frac{i\pi}{k^2-1}}\cot{\frac{ki\pi}{k^2-1}}&=&1+\cot{\left(\frac
{i\pi}{k^2-1}+\frac{ki\pi}{k^2-1}\right)}
\left(\cot{\frac{i\pi}{k^2-1}}+\cot{\frac{ki\pi}{k^2-1}}\right)\\
&=&1+\cot{\frac{i\pi}{k-1}}\left(\cot{\frac{i\pi}{k^2-1}}+\cot{\frac
{ki\pi}{k^2-1}}\right).
\end{eqnarray*}
For each $j\in\{1, 2, \cdots, k-2\}$, we regroup the summation by
if $i= j+(k-1)m,$ where $0\leq m\leq k$, it gives that
\begin{eqnarray*}
\sum_{m=0}^{k}\cot{\frac{i\pi}{k-1}}\cot{\frac{i\pi}{k^2-1}}&=&\cot
{\frac{j\pi}{k-1}}\left(\sum_{m=0}^{k}\cot{\frac{(j+(k-1)m)\pi}
{k^2-1}}\right)\\
&=&\cot{\frac{j\pi}{k-1}}\left(\sum_{m=0}^{k}\cot{\left(\frac{j\pi}
{k^2-1}+\frac{m\pi}{k+1}\right)}\right)\\
&=&(k+1)\cot^2{\frac{j\pi}{k-1}},
\end{eqnarray*}
where we use Lemma 1.2. by taking $x=\frac{j\pi}{k^2-1}$. And
similarly
\begin{eqnarray*}
\sum_{m=0}^{k}\cot{\frac{i\pi}{k-1}}\cot{\frac{ki\pi}{k^2-1}}&=&\cot
{\frac{j\pi}{k-1}}\left(\sum_{m=0}^{k}\cot{\frac{k(j+(k-1)m)\pi}
{k^2-1}}\right)\\
&=&\cot{\frac{j\pi}{k-1}}\left(-\sum_{m=0}^{k}\cot{\left(-\frac{kj\pi}
{k^2-1}+\frac{m\pi}{k+1}\right)}\right)\\
&=&(k+1)\cot^2{\frac{j\pi}{k-1}}.
\end{eqnarray*}
When $i=(k-1)j$, where $1\leq j\leq k$, it gives that
\[
\sum_{j=1}^k\cot{\frac{(k-1)j\pi}{k^2-1}}\cot{\frac{k(k-1)j\pi}
{k^2-1}}=-\sum_{j=1}^k\cot^2{\frac{j\pi}{k+1}}.
\]
Sum all terms up, it gives that
\begin{eqnarray*}
\sum_{i=1}^{k^2-2}\cot{\frac{i\pi}{k^2-1}}\cot{\frac{ki\pi}{k^2-1}}&=&
(k^2-2-k)+2(k+1)\sum_{j=1}^{k-2}
\cot^2{\frac{j\pi}{k-1}}-\sum_{j=1}^k\cot^2{\frac{j\pi}{k+1}}\\
&=&(k^2-2-k)+2(k+1)\frac{(k-2)(k-3)}{3}-\frac{k(k-1)}{3}\\
&=&\frac{2}{3}k^3-2k^2+2.
\end{eqnarray*}
\end{proof}
{\bf Acknowledgments}: The author wish to thank his advisor Chen
Xiuxiong to introduce him into the program of the existence of
extremal metrics on K\"ahler surface. The author is also grateful
to Hua Zheng, Jeff Viaclovsky,  C. LeBrun and Brian Weber for some helpful
discussions.

whe@math.wisc.edu\\
Department of Mathematics, University of Wisconsin-Madison,\\
Madison, Wisconsin, 53706, USA
\end{document}